\documentclass[12pt]
{amsart}
\usepackage{amsfonts,amsmath, amsthm,amssymb,indentfirst}

\usepackage[english,russian]{babel}

\renewcommand{\ge}{\geqslant}

\newcommand{\R}{\mathbb R}

\begin{document}

\title{Мотивированное изложение доказательства теоремы Тверберга}
\author{В. Ретинский, А. Рябичев, А. Скопенков}
\thanks{В. Ретинский: Высшая Школа Экономики.
\newline
А. Рябичев: Высшая Школа Экономики, Независимый Московский Университет.
\newline
А. Скопенков: \texttt{https://users.mccme.ru/skopenko}.
Московский Физико-Технический Институт, Независимый Московский Университет. 
Частично поддержан грантом РФФИ No. 19-01-00169.}
\date{}
\maketitle

\begin{abstract}
We present a motivated exposition of the proof of the following Tverberg Theorem: 
{\it For every integers $d,r$ any $(d+1)(r-1)+1$ points in $\R^d$ can be decomposed into $r$ groups
such that all the $r$ convex hulls of the groups have a common point.}
The proof is by well-known reduction to the B\'ar\'any Theorem.
However, our exposition is easier to grasp because additional constructions (of an embedding $\R^d\subset\R^{d+1}$, of vectors $\varphi_{j,i}$ and statement of the Bara\'ny Theorem) are not introduced in advance in a non-motivated way, but naturally appear in an attempt to construct the required decomposition.
This attempt is based on rewriting several equalities between vectors as one equality between vectors of higher dimension.
\end{abstract}

\newtheorem{theorem}{Теорема}

\bigskip
Мы приводим мотивированное изложение доказательства теоремы Тверберга (см. формулировку  ниже).
По сути наше изложение аналогично \cite[\S8.3]{Ma02} (см.\ также \cite{BO97}).
Однако оно проще для восприятия, поскольку дополнительные построения  (вложение $\R^d\subset\R^{d+1}$, векторы $\varphi_{j,i}$ и цветная теорема Каратеодори, см. ниже) не вводятся немотивированно заранее, а  естественно возникают при попытке построить нужное разбиение.
Эта попытка основана на записи нескольких равенств между векторами в виде одного равенства между векторами большей размерности.
 
{\bf Выпуклой оболочкой}  конечного набора точек $p_1,\ldots,p_n\in\R^d$ называется
$$\left<p_1,\ldots,p_n\right>:=\{\alpha_1p_1+\ldots+\alpha_np_n\ :\ \alpha_1,\ldots,\alpha_n\ge0,\ \alpha_1+\ldots+\alpha_n=1\}.$$

\begin{theorem}[Радон]\label{t:rad} Для любого $d$ любые $d+2$ точки в $\R^d$ можно разбить на два множества, выпуклые оболочки которых пересекаются.
\end{theorem}

\begin{theorem}[Тверберг]\label{t:tver} Для любых $d,r$ любые $(d+1)(r-1)+1$ точки в $\R^d$ можно разбить на $r$ множеств, выпуклые оболочки которых имеют общую точку.
\end{theorem}

Сперва продемонстрируем одну из идей на примере доказательства теоремы Радона (т.е. частного случая теоремы Тверберга для $r=2$).
Обозначим $[n]:=\{1,2,\ldots,n\}$.

\begin{proof}[Доказательство теоремы Радона]
Обозначим данные точки через $p_1, p_2, \ldots, p_{d+2}$.
Достаточно привести разбиение $[d+2]=A_1\sqcup A_2$ и действительные числа
$\alpha_1, \alpha_2, \ldots, \alpha_{d+2}\ge0$, для которых
$$\sum_{i\in A_1} \alpha_i p_i = \sum_{i\in A_2} \alpha_i p_i \quad\text{и}\quad
\sum_{i\in A_1}\alpha_i = \sum_{i\in A_2}\alpha_i, \eqno (1)$$
причем $\sum_{i\in A_1} \alpha_i=1$.
Вместо последнего равенства единице достаточно потребовать отличность от нуля.
Обозначим
$$p^+_i:=(p_i,1)\in \R^{d+1},\quad i\in[d+2],\quad\text{и}\quad
S_j^+:=\sum_{i\in A_j} \alpha_i p_i^+,\quad j=1,2.$$
Тогда два равенства (1) равносильны одному равенству $S_1^+=S_2^+$.

Точки $p_1, p_2, \ldots, p_{d+2}\in\R^{d+1}$ линейно зависимы.
Значит, существуют $\mu_1, \mu_2, \ldots, \mu_{d+2}$, не все равные нулю, для которых
$\sum_{i=1}^{d+2} \mu_i p_i^+ = 0$.
Перенесем в правую часть все слагаемые с отрицательными $\mu_i$.
Получим, что следующее разбиение и числа --- искомые:
$$A_1:=\left\{i\in[d+2]\ :\ \mu_i\ge0\right\},\quad
A_2:=\left\{i\in[d+2]\ :\ \mu_i<0\right\} \quad\text{и}\quad \alpha_i=|\mu_i|.$$
\end{proof}

\begin{proof}[Доказательство теоремы Тверберга]
Докажем теорему для $r=3$, доказательство общего случая аналогично.
Обозначим данные точки через $p_1, p_2, \ldots, p_{2d+3}$.
Достаточно привести разбиение $[2d+3]=A_1\sqcup A_2 \sqcup A_3$ и действительные числа
$\alpha_1, \alpha_2, \ldots, \alpha_{2d+3}\ge0$, для которых
$$\sum_{i\in A_1} \alpha_i p_i = \sum_{i\in A_2} \alpha_i p_i = \sum_{i\in A_3} \alpha_i p_i
\quad\text{и}\quad \sum_{i\in A_1}\alpha_i = \sum_{i\in A_2}\alpha_i = \sum_{i\in A_3}\alpha_i, \eqno (1')$$
причем $\sum_{i\in A_1} \alpha_i=1$.
Вместо последнего равенства единице достаточно потребовать равенство $\sum_{i=1}^{2d+3}\alpha_i=1$
(разделим на 3 все $\alpha_i$).
Обозначим
$$p^+_i:=(p_i,1)\in \R^{d+1}\quad\text{и}\quad S_j^+:=\sum_{i\in A_j} \alpha_i p_i^+.$$
Тогда равенства (1') равносильны следующим:
$$S_1^+=S_2^+=S_3^+\quad\Leftrightarrow\quad
\begin{cases}S_1^+=S_3^+ \\ S_2^+=S_3^+\end{cases}
\quad\Leftrightarrow\quad
\begin{cases}1\cdot S_1^++0\cdot S_2^+-S_3^+=0 \\ 0\cdot S_1^++1\cdot S_2^+-S_3^+=0\end{cases}
\quad\Leftrightarrow$$
$$\Leftrightarrow\quad(S_1^+,0)+(0,S_2^+)+(-S_3^+,-S_3^+)=0.$$
Обозначим
$$\varphi_{1,i}=(p^+_i,0)\in\R^{2d+2},\quad \varphi_{2,i}=(0, p^+_i)\in\R^{2d+2}\quad\text{и}\quad
\varphi_{3,i}=(-p^+_i, -p^+_i)\in\R^{2d+2}.$$
Тогда равенство $S_1^+=S_2^+=S_3^+$ равносильны следующему:
$$\sum_{i\in A_1}\alpha_i\varphi_{1,i} + \sum_{i\in A_2}\alpha_i\varphi_{2,i} + 
\sum_{i\in A_3}\alpha_i\varphi_{3,i}=0. \eqno(2)$$
Осталось применить следующую теорему.


\begin{theorem}[Барань; цветная теорема Каратеодори; \cite{Ba82}, {\cite[\S8.2]{Ma02}}] Пусть точка $0\in\R^n$ лежит в выпуклой оболочке каждого из $n+1$ конечных множеств  $M_1, M_2, \ldots, M_{n+1}\subset\R^n$.
Тогда существуют точки $m_i\in M_i$, для которых $0\in\left<m_1, m_2, \ldots, m_{n+1}\right>$.
\end{theorem}

Так как $0=\dfrac{\varphi_{1,i}+\varphi_{2,i}+\varphi_{3,i}}{3}\in\left<\varphi_{1,i},\varphi_{2,i},\varphi_{3,i}\right>$, то можно применить эту теорему к $n=2d+2$
и $M_i=\{\varphi_{1,i},\varphi_{2,i},\varphi_{3,i}\}$.
Получим точки $m_i\in M_i$, для которых $0\in\left<m_1, m_2, \ldots, m_{n+1}\right>$.
Для $j=1,2,3$ обозначим $A_j:=\left\{i\in[2d+3]\ :\ m_i=\varphi_{j,i}\right\}$.
Так как $0\in\left<m_1, m_2, \ldots, m_{n+1}\right>$, то существуют
$\alpha_1, \alpha_2, \ldots, \alpha_{2d+3}\ge0$, для которых выполнено равенство (2) и $\sum_{i=1}^{2d+3}\alpha_i=1$.
\end{proof}

\end{document}